\documentclass[12pt]{article}
\usepackage{amsmath,amssymb,amsthm,amsfonts,amscd}
\usepackage[numeric]{amsrefs}
\input colordvi

\swapnumbers

\newcommand\D{\Delta}
\newcommand\G{\Gamma}
\newcommand\f{\frac}

\newcommand{\Z}{{\mathbb{Z}}}
\newcommand{\R}{{\mathbb{R}}}

\renewcommand\Re{\text{Re~}}

\renewcommand\({\left(}
\renewcommand\){\right)}

\newcommand{\gobble}[1]{}
  \newcommand{\rangeref}[2]{%
    \ref{#1}--\afterassignment\gobble\fam 0\ref{#2}%
  }

\begin{document}

\title{A method for computing general automorphic forms on general groups}
\author{Stephen D. Miller\thanks{Partially supported
by NSF grant DMS-0601009.  The author wishes to thank Avner Ash, Andrew Booker, David Farmer, Wilfried Schmid, and Harold Stark for their helpful discussions.}\\
\tt{miller@math.rutgers.edu}}
\date{July 12, 2009}

\maketitle
\begin{abstract}
This article describes a general method for computing automorphic forms using Voronoi-type summation formulas. It gives a numerical example where the technique is successful in quickly finding a cusp form on $GL(3,\Z)\backslash GL(3,\R)$, albeit one whose existence was already known as a Langlands lift.

\vspace{.2cm}

\noindent MSC 2000 classification codes: 11F55, 11M41, 11Y99
\end{abstract}

The purpose of this note is to describe a method for computing general automorphic forms.  I have carried out only limited computational tests so far, and have not discovered any new automorphic forms using it.  However, the method does identify some lifted cusp forms on $GL(3)$, and until recently was the only general method to compute an automorphic form on a higher rank group.  It generalizes the methods of Hejhal (e.g.~\cite{hejhal}) and others for Maass forms on $GL(2)$, but does not require Whittaker or Bessel functions.  Aside from trying to compute some of the same objects, this method is otherwise unrelated to the cohomological methods developed by Ash and others which use geometric data to compute special types of automorphic forms.

I shall begin with a highbrow version of the method, and later explain its concrete manifestations and how they relate to existing methods.  Every cuspidal automorphic representation for $\G\backslash G$, where $G$ is a arbitrary reductive Lie group and $\G$ an arithmetically defined subgroup, has associated {\em automorphic distributions}  as in \cite{voronoi}.  These are  viewed as  objects in $C^{-\infty}(N)$, the distributions on the maximal unipotent subgroup $N$ of  $G$.  Furthermore, they satisfy invariance properties, in particular under $N_\G = \G\cap N$.   Fourier  expansions for
 $L^2(N_\G\backslash N)$ and $C^{-\infty}(N_\G\backslash N)$  can be computed  in an algorithmic way using tools from representation theory, for example using \cite{Brezin:1970} or Kirilov's method of co-adjoint nilpotent orbits.  The invariance under the rest of the group $\G$ then provides many nontrivial (distributional) identities between these Fourier components.  So far this has been carried out for $GL(3)$ and $GL(4)$ (see, for example, \cite{voronoi,glnvoronoi}) and $Sp(4)$ (\cite{chen}).  The full set of identities characterizes automorphy (a kind of converse theorem), and when these identities are integrated against test functions, they yield identities involving invariants of the automorphic representation (such as Hecke eigenvalues).  The main idea of this note is to use  these relations to numerically identify automorphic forms.

The integrated identities for $G=GL(3,\R)$ and $\G=GL(3,\Z)$ were worked out in detail in \cite{voronoi}, where they are the analogs of the classical Voronoi summation formulas for $GL(2)$.  These $GL(2)$ formulas have the form
\begin{equation}\label{voronoigl2}
    \sum_{n\neq 0} a_n \,e(-na/c)\,f(n) \ \ = \ \ |c|\sum_{n\neq 0}
    \f{a_n}{|n|}\,e(n\bar{a}/c)\,F(n/c^2)     \ \ \ \ \hbox{(see \cite[(1.12)]{voronoi}),}
\end{equation}
where $a/c$ is a fraction in lowest terms, $a\bar{a}\equiv 1\pmod c$, $a_n$ are the coefficients of the form, $e(u)=e^{2\pi i u}$, $f$ is a test function, and $F$ a type of transform of it (for details please see the statement in \cite{voronoi}).  When $f$ is chosen to be the appropriate $K$-Bessel function, the relation (\ref{voronoigl2}) amounts to the equality between the Fourier expansions of a Maass form evaluated at two $\G$-equivalent points in the upper half plane.  The generalized Voronoi formulas are similar equalities, but for  specific Fourier components on $N_\G\backslash N$, and more general vectors in the automorphic representation (e.g.~not just the spherical vector, which corresponds to the $K$-Bessel function).

 Such a relation, with the $K$-Bessel function, has been the starting point for many numerical computations of Maass forms, which seek to solve identities such as (\ref{voronoigl2}) for the values of the coefficients $a_n$.  The essential point of the method here is that the Voronoi formula still affords this opportunity, and in fact with much greater flexibility in that it does not force a single choice of test function. Indeed, $f$ can be a certain fractional power times any Schwartz function.

 More importantly, Voronoi formulas have now been established for $GL(n)$ in \cite{voronoi,glnvoronoi}.
The $GL(3)$ formula has the following form:
\begin{equation}\label{voronoigl3}
     \sum_{n\neq 0} a_{q,n} \,e(-na/c)\,f(n) \ \ = \ \ |c|\sum_{d|cq}\sum_{n\neq 0} \f{a_{n,d}}{|nd|}\,S(q\bar{a},n;qc/d)\,F\(\f{nd^2}{c^3q}\),
\end{equation}
where $S(a,b;c)=\sum_{x\in(\Z/c\Z)^*}e(\f{ax+b\bar{x}}{c})$ is the usual Kloosterman sum, and $F$ represents a different transform of $f$ than in (\ref{voronoigl2}).  The Fourier coefficients on $GL(3)$, as they are on other higher rank groups, are indexed by more than one integral parameter, though the $L$-function is determined by the data $\{a_{1,n}=\overline{a_{n,1}}\}$.

  Thus one can simply seek to apply the existing coefficient-solving methods to formula (\ref{voronoigl3}).  This seems to be advantageous compared to long-standing proposals to explicitly compute $GL(n)$ Whittaker functions:~for one thing, it is uniform over all representations, and again subsumes the Whittaker approach.  Thus it cannot be any less successful, though choosing appropriate functions may still be a challenge.  Also, the Fourier expansion in higher rank involves more complicated sums than it does for classical modular or Maass forms, and so it is appealing to isolate individual Fourier components of it like distributional identities such as (\ref{voronoigl3}) do.
  The coefficient-solving methods also include techniques to hone in on the correct eigenvalue of form when it is not known in advance, which in principal apply here also.  I have not been able to investigate this in detail, however.

 I have attempted to carry out the method  in a simple case, of {\em cohomological} cusp forms on $GL(3,\Z)\backslash GL(3,\R)$.  The reason for doing so is that the representation type is pinned down exactly, and the coefficients are also expected to have integrality properties.  Indeed, for many  purposes cohomological forms are the appropriate analog of discrete series/holomorphic forms on $GL(2)$. I wish to stress that this is done simply to avoid precision problems, and that the technique is in principal general.  The simplest case of (\ref{voronoigl3}) occurs when there is no additive twist, corresponding to $a=0$ and $c=1$:
 \begin{equation}\label{simplevoronoi}
    \sum a_n \,f(n) \ \ = \ \ \sum \overline{a_n}\, F(n)\,,
 \end{equation}
 where again $f$ and $F$ are a transform pair of functions.  This formula can also be proven by writing the lefthand side as a contour integral in terms of the Mellin transform of $f$,
  $\f{1}{2\pi i}\int_{\Re{s}=2}L(s)Mf(s)ds$, and applying the functional equation of the $L$-function $L(s)=\sum a_n n^{-s}$.  Formula (\ref{simplevoronoi}) is thus essentially an equivalent form of the approximate functional equation.  One could also of course attempt to use the additive twists $e(-na/c)$ present in (\ref{voronoigl3}) to gain further information.
    In the experiments $f(x)$ was chosen to be a polynomial times a Gaussian, and $F$ was computed using \cite[(5.27)]{voronoi} -- more accurately, approximated by a partial sum of its residues since a closed form expression for $F$ was impractical.  An important difference between higher and lower rank groups is that more coefficients are needed as $c$ grows, owing to the $c^d$ growth of the twisted conductor of a degree $d$ $L$-function, such as one coming from $GL(d)$.  The smaller growth for $GL(2)$ is an important ingredient in Stark's Hecke operator method \cite{stark} (see also \cite[pp.220-224]{terras}), which accordingly appears difficult to generalize to higher rank.

  The experiments were carried out to search for cusp forms whose representation type agrees with the symmetric square lift of Ramanujan's $\D$ cusp form of weight 12.  Here one knows via earlier computations (originating in studying the group cohomology of $SL(3,\Z)$ -- see \cite{ashpol}, for example) that there is only one such form to be found:~the global lift of $\Delta$.  Locking in the value $a_1=1$ and using $a_2,\ldots,a_{50}$ as 49 unknowns, 49 different Voronoi equations were created.  The truncation of identity (\ref{simplevoronoi})  to its first 50 terms on each side  was then solved for $a_2,\ldots,a_{50}$ using matrix inversion.  The first few coefficients found this way are listed in Table~\ref{table1}, along with their correct values  (which can be independently calculated using the description of the symmetric square lift).
  \begin{table}\centering
  \begin{tabular}{|r||c|c|}
    \hline
    $n$ &  Approximate Value & Exact Value \\[3pt]
\hline
    2 & -0.718812 & $-\f{23}{32}\approx$ -0.718750 \\[3pt]
    3 & -0.641492  & $-\f{1403}{2187} \approx$ -0.641518 \\[3pt]
    4 &  1.235340  & $\f{1265}{1024}\approx$ 1.23535 \\[3pt]
    5 & -0.522150  & $-\f{1019969}{1953123}\approx$ -0.522224 \\[3pt]
    6 &  0.460477  & $\f{32269}{69984}\approx$ 0.461091 \\[3pt]
    7 & -0.854928  & $ -\f{34631943}{40353607}\approx$ -0.858212 \\[3pt]
    8 & -0.418517  & $-\f{13255}{32768}\approx$ -0.40451 \\[3pt]
    \hline
  \end{tabular}\caption{Coefficients of a cohomological cusp form on $GL(3,\Z)\backslash GL(3,\R)$.  This form is the symmetric square of Ramanujan's $\Delta$ cusp form.  The computation took 8 minutes on one core of an Intel(R) Core 2 CPU system running at 1.5 GHz, with 2 GB of RAM, slower than the methods of \cite{ashpol}.}\label{table1}
  \end{table}
  One could further use the
  fact that the correct values of $a_n n^{11}$ are integers to eliminate precision errors (in this atypical situation).

  It is natural to ask whether or not there are cohomological cusp forms for $GL(3,\Z)\backslash GL(3,\R)$ aside from these symmetric square lifts from $GL(2)$.  They do exist for general congruence subgroups, but none are known at full level -- this despite a serious numerical search by Ash and Pollack \cite{ashpol}, who consequently conjectured that none exist.  (See also \cite{calmaz}, who propose another explanation for the paucity of these forms.)

 In the case of Maass-type (i.e.~spherical) cusp forms on $GL(3)$, one does not know the eigenvalue to start with, and thus must apply procedures similar to those developed for the $GL(2)$ case to find it along with the coefficients.  Though there are again symmetric square lifts of $GL(2)$ Maass forms, it was proven in \cite{thesis} that there do exist many nonlifted forms:~in fact the lifted ones comprise only a minuscule portion of the spherical automorphic spectrum on $GL(n)$ \cite{muller}.  Nevertheless, none have been identified to date, and finding them remains an important challenge.

\subsection*{Recent Developments}

Shortly after this note was first distributed on the arXiv, significant progress was made in computing Maass-type   automorphic forms on $GL(3)$.  The methods used are similar to those described above in that they use an integrated form of the functional equation, as opposed to trying to use the Whittaker expansion directly.  However, the success of these methods depends very much on some important differences.

 First, Andrew Booker and Ce Bian announced the breakthrough calculation of new automorphic forms for $SL(3,\Z)\backslash SL(3,\R)/SO(3,\R)$ \cite{booker}.  Their method uses the approximate functional equation (\ref{simplevoronoi}), but also for the $L$-function twisted by Dirichlet characters.  This is equivalent to using the additively-twisted Voronoi formula (\ref{voronoigl3}), though it is not clear yet whether grouping into additive or multiplicative twists is more efficient computationally.  A key point in their work is creating sufficient numerical stability for their large system of roughly 10,000 equations (related to \cite{bookerthesis}).  The payoff is clear:~their results go well beyond the experiment above, which used only 49 equations.   The use of twisted functional equations is motivated by the $GL(3)\times GL(1)$ converse theorem of \cite{jpss}, which shows that these $L$-functions are sufficient to characterize automorphy.    The identity (\ref{voronoigl3}) shares this feature, as do the generalizations to $C^{-\infty}(N_\G\backslash G)$ Fourier expansions mentioned earlier.  For this reason, it may be possible to use the distributional identities in \cite{chen} to compute cusp forms for $Sp(4)$, e.g.~Siegel modular forms.

Days later, David Farmer, Sally Koutsoliotas, and Stefan Lemurell announced another method to numerically compute automorphic forms on
$SL(3,\Z)\backslash SL(3,\R)/SO(3,\R)$ \cite{farmer}.  They again use the approximate functional equation (\ref{simplevoronoi}), with a wider choice of test functions.  However, a crucial point is that they force certain expected Hecke relations (such as $a_2a_3=a_6$ and $a_2a_5=a_{10}$), thereby obtaining a nonlinear system of equations.  They report that this has the effect of drastically lowering the number of equations needed down to about 20 for the smallest examples.  Despite this reduced number of equations, their numerical accuracy still exceeds that of Table~\ref{table1}.

Finally, a third method was introduced by Borislav Mezhericher \cite{boris}.
This method develops a formula which is simultaneously both a consequence of (\ref{voronoigl3}), and yet a generalization in that it reduces to (\ref{voronoigl3}) when an auxiliary parameter is set to zero.  It shares some features in common with Hejhal's $GL(2)$ approach, and may be useful for employing techniques that are successful there.

It should be noted that, in computing spherical automorphic forms for $GL(n,\Z)$, there are no integrality results that can be relied upon:~every parameter of these forms, whether they be coefficients or eigenvalues, is thought to be transcendental.  Thus all methods include an approach to identify the correct eigenvalue parameters.

\end{document}